\documentclass[11pt]{amsart}

\usepackage{a4,amssymb,amsmath,amsthm,xspace,epsfig,latexsym}

\newcommand{\Chat}{\ensuremath{\widehat{\mathbb{C}}}}

\newtheorem{thm}{Theorem}
\newtheorem{lemma}{Lemma}
\newtheorem{corollary}{Corollary}
\newtheorem*{theorem*}{Theorem}

\newenvironment{dem}{\begin{proof}[{\sc Proof}]}{\end{proof}}

\title{Mean value conjectures for rational maps}

\author{Edward Crane\\ Merton College, Oxford}

\begin{document}

\begin{abstract} 
Let $p$ be a polynomial in one complex variable. Smale's mean value conjecture estimates $|p'(z)|$ in terms of the gradient of a chord from $(z, p(z))$ to some stationary point on the graph of $p$. The conjecture does not immediately generalise to rational maps since its formulation is invariant under the group of affine maps, not the full M\"{o}bius group. Here we give two possible generalisations to rational maps, both of which are M\"{o}bius invariant. In both cases we prove a version with a weaker constant, in parallel to the situation for Smale's mean value conjecture. Finally, we discuss some candidate extremal rational maps, namely rational maps all of whose critical points are fixed points.
\end{abstract}
\maketitle
\section{Introduction}

 Let $p$ be any polynomial with coefficients in $\mathbb{C}$. Then $\zeta \in \mathbb{C}$ is a \emph{critical point} of $p$ when $p'(\zeta) = 0$. Its image $p(\zeta)$ is the corresponding \emph{critical value}. In 1981 Steve Smale proved the following result about critical points and critical values of polynomials, in connection with algorithms for finding roots of polynomials.
\begin{theorem*}[Smale, \cite{Sm}]\label{T: Smale}{\quad}\\ Let $p$ be a polynomial of degree $n \ge 2$ over $\mathbb{C}$ and suppose that $x \in \mathbb{C}$ is not a critical point of $p$. Then there exists a critical point $\zeta$ of $p$ such that \[\left| \frac{p(\zeta) - p(x)}{\zeta - x}\right| \,\le\, 4\, |p'(x)|\, .\]
\end{theorem*}
Thus the derivative of $p$ at each point can be estimated in terms of the gradients of chords to the critical points; in this way it is similar to the ordinary mean value theorem. Letting $\zeta_1, \dots, \zeta_{n-1}$ be the critical points of $p$, repeated according to multiplicity. We define \[S(p,x) = \min_i \left|\frac{p(\zeta_i) - p(x)}{(\zeta_i - x)p'(x)}\right|\] Smale conjectured that the correct bound for polynomials of degree $n$ is \begin{equation}\label{E: SMVC} S(p,x) \le 1 - 1/n\,,\end{equation} which would be best possible in view of the example \[ p_0(z) = (z^n - nz)/(1-n)\, \, \text{ with }\, x=0\,.\] 
 It is known to be the correct bound for various special classes of polynomials, and when $n=2,3,4,5$. See \cite{Crane},\cite{RaSch}, \cite{Sch}, \cite{Sh-S}, \cite{Tis1}, and \cite{Tis2} for details. The best known upper bound for general polynomials of general degree $n$ is $S(p,x) \le 4^{1 - 1/(n-1)}$, which is established in \cite{BMN}.

 Smale's mean value conjecture does not immediately generalise to rational maps because the point at infinity in the Riemann sphere, $\infty \in \Chat$, is treated in a special way. Indeed the quantity $S(p,x)$ is invariant under the action of the group of affine maps $z \mapsto az + b$, which is the stabiliser of $\infty$ in the action of the M\"obius group on $\Chat$. If $A$ and $B$ are two affine maps then $S\left(A \circ p \circ B, B^{-1}x\right) = S(p,x)$. Any reasonable generalisation to rational maps should have a formulation that is invariant under the full M\"obius group. In \S \ref{S: direct generalisation} we give such a generalisation, and in \S \ref{S: symmetric generalisation} we prove a similar result which looks more natural but does not include Smale's theorem as a special case. 
 
 In \S \ref{S: rational maps with critical points fixed}, we discuss some candidate extrema, namely rational maps all of whose critical points are fixed points. After discussing their construction and connection with Newton's method, we investigate the possible values of the multipliers at their non-critical fixed points, using hyperbolic geometry and the rational fixed point index theorem.
 
 The results in this paper appear in the author's PhD thesis \cite{thesis}. 

\section{A direct generalisation}\label{S: direct generalisation}

The proof of Smale's theorem only uses one fact about polynomials that does not apply to general rational maps, namely that $P(\infty)=\infty$. It does not use the facts that that $\infty$ is a critical point for $P$ or that there are no other pre-images of $\infty$. This allows us to use the same method to prove a generalisation for rational maps. For a rational map $R$ and any point $x \in \Chat$ such that $R(x)=x$, we denote by $R^\#(x)$ the derivative of $R$ at $x$ with respect to any local co-ordinate at $x$; this is also called the \emph{multiplier} of $R$ at $x$.

\begin{thm}\label{T: one critical point}{\quad}\\ Let $R$ be any rational map of degree at least 2, and let $x$, $y$ be points of $\Chat$ with $R(x) \neq R(y)$, such that $x$ is not a critical point of $R$. Then there exists a critical point $\zeta$ of $R$ and a M\"obius map $M$ such that $M \circ R$ fixes each of $ x$, $y$ and $\zeta$, and  $|(M \circ R)^\#(x)| > 1/4$.
\end{thm}

\begin{dem} The statement is true for $R$ if and only if it is true for $\hat{R} = T \circ R \circ S$, where $T$ and $S$ are any M\"obius maps. We choose a M\"obius map $S$ such that $S(\infty) = y$ and $S(0) = x$, and a M\"obius map $T$ such that $T(R(y)) = \infty$ and $T(R(x)) = 0$. Then $\hat{R} = T \circ R \circ S$ is a rational map fixing $0$ and $\infty$. Let $D$ be the largest open disc centred on $0$ that carries a single-valued branch $\beta$ of $\hat{R}^{-1}$ that maps $0$ to $0$. Let $E = \beta(D)$; this is the component of $\hat{R}^{-1}(D)$ containing $0$. There is a critical point $\zeta$ of $\hat{R}$ on the boundary $\partial E$, and $\hat{R}(\zeta) \in \partial D$. Now $\beta$ is a univalent map from $D(0, |\hat{R}(\zeta)|)$ to the domain $E$ and it omits the value $\zeta$. So Koebe's $\frac{1}{4}$-Theorem yields \[|(\hat{R}^{-1})'(0)| \,\le\, 4 \,\frac{|\zeta|}{|\hat{R}(\zeta)|}\, .\]  
The equality case of Koebe's theorem cannot occur, because then the branch $\beta$ of $R^{-1}$ would cover the whole of the Riemann sphere except for a slit, and the complement of $D$ could not be contained in the image of $R$.
\end{dem}

 The special case with $y=\infty$ and $R$ a polynomial is a simple reformulation of Smale's result. An alternative formulation of the conclusion, valid when none of the relevant points is $\infty$, is that there exists a critical point $\zeta$ such that 
\[ \left|\frac{R'(x) (y-x)}{R(y)-R(x)}\right|\, > \, 
\frac{1}{4} \left|\frac{R(x)-R(\zeta)}{R(y)-R(\zeta)} \cdot \frac{y - \zeta}{x-\zeta}\right|  \, = \, \frac{1}{4} \left|\frac{R(x)-R(\zeta)}{x-\zeta} \cdot \frac{y - \zeta}{R(y)-R(\zeta)}\right| \,.       \]
For a fixed choice of $R$, when $x$ and $y$ are close together and not close to any critical points or poles of $R$, one would expect the left-hand side of this inequality to be close to $1$, and the middle expression should be close to $1/4$ so long as we choose a critical point $\zeta$ such that $R(\zeta)$ is not accidentally close to $R(x)$ or $R(y)$. On the other hand if the left-hand side is very small, this can be explained by the fact that $x$ is close to a critical point, and the first factor in the rightmost expression will be small. In this formulation, Theorem \ref{T: one critical point} certainly deserves to be called a mean value theorem.

 We ask what is the largest possible constant $K_1$ that can replace $1/4$ in Theorem \ref{T: one critical point}, possibly in terms of $n=\deg(R)$. Consider the example $R(z) = p_0(z) = (z^n - nz)/(1-n)$, $x=0$, $y = \infty$. Since $p(0)=0$, $p(\infty)=\infty$ and each finite critical point is fixed, the M\"obius map $M$ in Theorem \ref{T: one critical point} must be the identity regardless of the choice of $\zeta$. This shows that $K_1 \le n/(n-1)$. We conjecture that $K_1=n/(n-1)$. 

A rational map of degree $2$ must have exactly two simple critical points; after composing with M\"obius maps  we may assume that the map is $z \mapsto z^2$, and that $x = 1$ in the theorem; in this case the derivative at $x$ is equal to $2$. So the situation for degree $2$ is no different from the polynomial case. For degree $3$, the parameter space for the problem is 3-dimensional over $\mathbb{C}$, so it is unlikely that there is an elementary treatment, but a computer proof may be feasible.

\section{A symmetric generalisation} \label{S: symmetric generalisation}

Theorem \ref{T: one critical point} deals with the location of one critical point and critical value relative to two arbitrarily given points. On the other Smale's mean value conjecture deals with the location of two critical points and critical values relative to one arbitrary point point. We now prove a theorem that generalises this aspect to the setting of rational maps, though it does not contain Smale's Theorem as a special case because it no longer makes sense to specify that one of the critical points be $\infty$.
\pagebreak
\begin{thm}\label{T: two critical points}{\quad}\\
Let $R$ be a rational map of degree at least $2$, and $x \in \Chat$ be given, such that $R(x)$ is not a critical value of $R$. Then there exist two critical points $\zeta$ and $\kappa$ of $R$ such that if $M$ is the M\"obius map that makes $x$, $\zeta$ and $\kappa$ fixed points of $M \circ R$, then $|(M \circ R)^{\#}(x)| \ge 1/2$.
\end{thm}
\begin{dem}
The conclusion is unchanged if we replace $R$ by $S \circ R \circ T$, where $S$ and $T$ are any M\"obius maps, so w.\ l.\ o.\ g.\  we may assume $x = \infty = R(x)$; then all the critical values of $R$ are finite and there are at least two of them. Let $K$ be the convex hull of the critical values. Then there exist two critical points $\zeta$ and $\kappa$ such that $\text{diam}(K) = |R(\zeta) - R(\kappa)| > 0$.   Now there exists a single-valued branch $\beta$ of $R^{-1}$ defined on $\Chat \setminus K$, taking $\infty$ to $\infty$, and omitting $\zeta$ and $\kappa$. Let $f: \Chat \setminus \overline{\mathbb{D}} \to \Chat \setminus K$ be a Riemann map fixing $\infty$. The logarithmic capacity of a compact subset of $\mathbb{C}$ of diameter $d$ is at most $d/2$ (see \cite[Theorem 5.3.4]{Ran}). Therefore $f^{\#}(\infty) \ge 2 / |R(\zeta) - R(\kappa)|$. The Koebe $\frac{1}{4}$-Theorem tells us that $(\beta \circ f)^{\#}(\infty) \le 4 /|\zeta - \kappa|$, so we have $\frac{1}{R^{\#}(\infty)} = \beta^{\#}(\infty) \le 2 \frac{|R(\zeta) - R(\kappa)|}{|\zeta -\kappa|}$. Take $M$ to be the affine map that makes $M \circ R$ fix $\zeta$, $\kappa$ and $\infty$. Then $M^{\#}(\infty) = \frac{|R(\zeta) - R(\kappa)|}{|\zeta - \kappa|}$, which gives the result. 
\end{dem}

Now we ask what is the largest possible constant $K_2$ that can replace $1/2$ in Theorem \ref{T: two critical points}. The same example as for Theorem \ref{T: one critical point} shows that $K_2 \le n/(n-1)$ and we conjecture that $K_2=n/(n-1)$. 

\section{Rational maps with all critical points fixed}\label{S: rational maps with critical points fixed}

\subsection{A strategy for bounding $K_1$ and $K_2$ above}

 If $R$ is any rational map for which every critical point is a fixed point, then the absolute value of the multiplier at any other fixed point provides an upper bound for $K_1$ and $K_2$, since the M\"obius map $M$ appearing in Theorems \ref{T: one critical point} and \ref{T: two critical points} must be the identity. 
 
  A rational map of degree $n \ge 2$ has $2n-2$ critical points and $n+1$ fixed points, both counted with multiplicity. In this section we give a construction for rational maps all of whose critical points are fixed, and which have precisely one further fixed point. So in general they will have multiple critical points. Note that there also exist rational maps all of whose critical points are fixed but which have more than one additional fixed point, for example $z \mapsto z^n$ for $n \ge 3$, (which only gives the bounds $K_1, K_2 \le n$).

Let $g$ be any rational function. Then we define $R_g$ to be the associated Newton-Raphson map
\[R_g(z)\, = \,z - \frac{g(z)}{g'(z)}\,.\]
$R_g$ is of course designed so that it has a fixed point at each of the roots and poles of $g$; this fixed point is attracting in the case of a root and superattracting for a simple root. So  the roots of $g$ can be approximated by iterating $R_g$ from suitable initial values. In fact \[R_g'(z) = \frac{g(z) g''(z)}{g'(z)^2}\,.\] 
Suppose we can find a polynomial $h$ of degree $n$ such that $h$ has no repeated roots and all the roots of $h''$ are also roots of $h$. Then $h'$ has no multiple roots, so $R_h$ has no multiple poles. In particular, \emph{every critical point of $R_h$ is a fixed point}. Note that $R_h$ has a fixed point at $\infty$ with multiplier $R_h^{\#}(\infty) = \frac{n}{n-1}$. So such an example gives an alternative proof that $K_1, K_2 \le n/(n-1)$. For example, we could take $h(z) = a(z-c)^n + b(z-c)$, and then we find that $R_h$ is M\"obius-conjugate to the polynomial $p_0$, (via $z \mapsto 1/z$), so we recover the original example this way. It turns out that this is the only possibility for $h$ such that $R_h$ is conjugate to a polynomial.

For another example, take \[h(z) = z^4 + 2 z^3 + 6 z^2 + 5 z + 4 = (z^2 + z + 4)(z^2 + z + 1)\,.\]  In this case we have 
\[R_h(z) = \frac{3z^4 + 4z^3 + 6z^2 - 4}{4z^3 + 6z^2 + 12z + 5}\,\] \[h''(z) = 12(z^2 + z + 1)\,.\] Since $h$ has no repeated roots, the roots of $h''$ are roots of $h$ but not of $h'$. Thus $R_h$ has two finite fixed points of valency three, two finite fixed points of valency two, and no other critical points. $R_h$ is not M\"obius-equivalent to a polynomial because it does not have a critical point of valency $4$.

So in Theorems \ref{T: one critical point} and \ref{T: two critical points}, it seems likely that there are several extremal rational maps for each degree, even up to M\"obius equivalence. In contrast, for Smale's mean value conjecture it is conjectured that the extremal polynomial is unique up to affine equivalence, and this is known for $n\le 5$. Any method of proof for Smale's mean value conjecture that might also deal with our rational map conjectures ought not to rely on the uniqueness of the extremum.

\subsection{The Rational Fixed Point Theorem}

We will use the notion of residue fixed point index for holomorphic maps. If $f: U \to \mathbb{C}$ is holomorphic on an open set $U \subset \mathbb{C}$ and $z_0$ is an isolated fixed point of $f$, then the residue fixed point index of $f$ at $z_0$ is defined as \[\iota(f,z_0) = \frac{1}{2 \pi i} \int \frac{dz}{z - f(z)} \; ,\]
where the integral is taken around a positively-oriented circle around $z_0$ so small that it contains no other fixed points of $f$. If $z_0$ is a simple fixed point of $f$ then \[\iota(f,z_0) = \frac{1}{1-f^{\#}(z_0)}\,.\] If it is not simple then the index is still well-defined and finite, but this formula does not apply. The following theorem is an easy consequence of Cauchy's residue theorem.
\begin{theorem*}[Rational Fixed Point Theorem]{\quad} \\ For any rational map $f: \Chat \to \Chat$ which is not the identity map, the sum over all the fixed points of the residue fixed point index is $1$.\end{theorem*}
 See \cite[Chapter 12]{Mil} for a full discussion of this material, with applications to complex dynamics.
 
\subsection{Characterisation of Newton-Raphson maps}

The following lemma characterises those rational maps that occur as Newton-Raphson maps of rational functions in terms of their fixed points and the multipliers at those fixed points. Buff and Henriksen \cite{BuHe} characterised the rational maps that occur as K\"onig's methods for polynomials; this includes Newton's method for polynomials as a special case for which they cite \cite[Prop. 2.1.2]{He}. The extension here to Newton maps of rational functions may be new.  
\begin{lemma}[Characterisation of Newton-Raphson maps of rational functions]{\quad}\\
A rational map $R$ is the Newton-Raphson map associated to some non-linear rational function $g$ if and only if all the fixed points of $R$ are simple and each fixed point $\chi$ of $R$ satisfies $\frac{1}{1 - R^{\#}(\chi)} \in \mathbb{Z}$. Moreover, $g$ is a polynomial if and only if these integers are all positive.
\end{lemma}
A fixed point is called simple if its multiplier is not equal to $1$. It is a consequence of the rational fixed point theorem that one need only check the integrality condition for all but one fixed point. $R$ has a simple non-critical fixed point at $\infty$ precisely when \[\lim_{z \to \infty} \frac{R(z)}{z} \in \Chat \setminus \{0,1, \infty\}.\]
\begin{dem}
Suppose that $g(z) = p(z)/q(z)$ in lowest terms is a rational function of degree $n \ge 2$ and $R = R_g(z):= z - \frac{g(z)}{g'(z)}$. The fixed points of $R$ are precisely the roots and poles $\chi$ of $g$. If $g$ has order $m \neq 0$ at $\chi \in \mathbb{C}$ then \[1 - R'(\chi) = \left.{\frac{\textup{d}}{\textup{d}z}\left(\frac{g(z)}{g'(z)} \right)}\right|_{z = \chi} = \frac{1}{m}\,.\]
 For the converse, consider \[g(z):= \exp \int^z \frac{\textup{d}z}{z - R(z)}\,.\] This function is certainly locally defined away from fixed points of $R$ and satisfies the Newton-Raphson equation $R(z) = z - \frac{g(z)}{g'(z)}$. The integrand $\frac{1}{z - R(z)}$ has no multiple poles since at each fixed point of $R$ we are told $\frac{\textup{d}}{\textup{d}z}(z - R(z)) \neq 0$. We can therefore express the integrand in partial fractions as \[\frac{1}{z - R(z)} = q(z) + \sum_i \frac{A_i}{z-\chi_i} \,,\] where $q$ is a polynomial. $R$ does not have a multiple fixed point at $\infty$, so $z - R(z) \to \infty$ as $z \to \infty$, hence $q = 0$. Near a fixed point $\chi_i$ of $R$ we know that $z - R(z) = (z-\chi_i)/m + O((z-\chi_i)^2)$ for some positive integer $m_i$, so $A_i = m_i$. Now we can perform the integration explicitly: \[ \int^z \frac{\textup{d}z}{z - R(z)} = \sum_i m_i \log(z-\chi_i) + c\,,\] so \[g(z) = \exp(c). \prod_i (z-\chi_i)^{m_i},\] which is a rational function; it is a polynomial when all $m_i \ge 1$. 
\end{dem}

The following special case is also a special case of \cite[Prop. 4]{BuHe}.
\begin{corollary}[Characterisation of Newton-Raphson maps associated to polynomials without repeated roots]\label{C: no repeated roots}{\quad}\\
The following are equivalent for a rational map $R$:
\begin{enumerate}
\item $R$ has a simple fixed point at $\infty$ and all the finite fixed points of $R$ are also critical points of $R$;
\item $R(z) = z - \frac{g(z)}{g'(z)}$ for some non-linear polynomial $g$ with no repeated roots.
\end{enumerate}
In this situation, the fixed point of $R$ at $\infty$ is not critical; indeed its multiplier is $n/(n-1)$, where $n = \deg g = \deg R$.
\end{corollary} 

Suppose that a Newton-Raphson map $R_h$ associated to a polynomial $h$ has all its critical points fixed. Then we claim that all of its finite fixed points must be critical. Indeed, $R_h$ would otherwise have a non-critical fixed point of multiplier less than $1$ in modulus. The iterates of $R_h$ have exactly the same critical points as $R_h$, so some iterate of $R_h$ would violate Theorem \ref{T: two critical points}, which proves the claim. Corollary \ref{C: no repeated roots} then applies to show that $h$ has no repeated roots, and so every root of $h''$ is a critical point of $R_h$, hence a fixed point of $R_h$, so a root of $h$. Thus our original conditions on the polynomial $h$ are necessary as well as sufficient to ensure that all the critical points of $R_h$ are fixed.

For the Newton-Raphson map $R_g$ associated to a \emph{rational} function $g$ of degree $n$, we could have repelling fixed points associated to poles of $g$, but the multiplier would be $k/(k-1)$ for a pole of $g$ of order $k$. It is possible for the degree of $R_g$ to be less than the degree of $g$, so one might still hope to produce counterexamples to our conjectures using the present construction; we rule this out in the next section. 

\subsection{Forbidden multipliers}

The special case of Smale's mean value conjecture for polynomials in which the critical points are all fixed has been called Kostrikin's conjecture. Shub observed that if the critical points of a polynomial are all fixed, then the multiplier at each remaining fixed point must be greater than or equal to $1$ in modulus, and we used the same iterative argument above. We now give a different argument which also applies to rational maps whose critical points are all fixed. It shows that the multiplier must be strictly greater than $1$ in modulus.

 Consider first the case in which $R$ has only two critical values. Then $R$ pulls back a complete Euclidean metric from the conformal cylinder $\Chat \setminus \{ \text{critical values}\}$ to $\Chat \setminus R^{-1}(\{\text{critical values}\})$, so there can only be two critical points, and the map is M\"obius-conjugate to $z \mapsto z^n$, in which the multiplier at the remaining fixed points is $n$. 
 
 Now suppose that there are at least three critical values. Then any branch of the algebraic function $R^{-1}: \Chat \setminus \{ \text{critical values}\} \to \Chat \setminus R^{-1}(\{\text{critical values}\})$ can be lifted and continued to give a conformal isomorphism between the universal covers of these two domains. So any branch of $R^{-1}$ is locally a hyperbolic isometry between the complete hyperbolic metrics on $\Chat \setminus \{ \text{critical values}\}$ and $\Chat \setminus R^{-1}(\{\text{critical values}\})$. However, the inclusion map $I: \Chat \setminus \{ \text{critical values}\} \to \Chat \setminus R^{-1}(\{\text{critical values}\})$ omits some points because it is impossible for each of the critical values to have only one pre-image. Therefore the inclusion map is everywhere a strict contraction between the hyperbolic metrics. At any non-critical fixed point this shows that the multiplier is greater than 1 in modulus.

The following theorem excludes further values of the multiplier at any non-critical fixed point, so it makes a small amount of progress on Kostrikin's Conjecture.

\begin{thm}\label{T: excluded multipliers}{\quad}\\ 
Suppose that $R$ is a rational map all of whose critical points are fixed. Suppose that $R$ has degree $n$ and has $m$ critical points (not counting multiplicity). Then $R$ has exactly $n+1$ fixed points, of which $n+1-m$ are non-critical. The multiplier at any non-critical point does not lie in the closed disc whose diameter is the interval $[1, 1+ \frac{2}{n+m-2}]$, except in the case where $m=n$ and the multiplier of the remaining fixed point is $n/(n-1)$, which is on the boundary of this disc.
\end{thm}
\begin{dem}
 We showed above that any non-critical fixed points must be repelling, hence simple, so there are no multiple fixed points of $R$, and hence there are exactly $n+1$ distinct fixed points. The residue fixed point index of a fixed critical point is $1$.  In the case where there are $n$ critical fixed points, the rational fixed point theorem shows that the remaining fixed point must have residue fixed point index $1-n$, so must have multiplier $n/(n-1)$. In the general case, note that for a non-critical fixed point, the multiplier has modulus strictly greater than one, so its residue fixed point index has real part strictly less than 1/2. Select a particular non-critical fixed point $z_0$. We get a contribution of less than $m + \frac{(n-m)}{2}$ to the real part of the total index from all the other fixed points, so the residue fixed point index $\iota(R,z_0)$ has real part greater than $1 - (m + \frac{(n-m)}{2}) = 1 - \frac{m+n}{2}$. Hence the multiplier does not lie in the disc whose diameter is the interval $[1, 1 + \frac{2}{m+n-2}]$, as required. In fact, further small regions of excluded values of the multiplier may be found by considering iterates of $R$.
\end{dem} 
 In the case $m < n$, this disc of excluded values contains the multiplier $\frac{n}{n-1}$ in its interior.  Since the multipliers at the fixed points of a Newton-Raphson map associated to a rational function are real and non-negative, Theorem \ref{T: excluded multipliers} implies that one cannot improve on the multiplier $n/(n-1)$ by using Newton-Raphson maps associated to rational functions. Any counterexample for our conjectures will have to come from another source.

\end{document}